\documentstyle[12pt]{article}
\newcommand{\const}{\mathop{\rm const}\limits}

\newcommand{\supp}{\mathop{\rm supp}\limits}

\newcommand{\Var}{\mathop{\rm Var}\limits}

\newcommand{\rank}{\mathop{\rm rank}\limits}

\newcommand{\diam}{\mathop{\rm diam}\limits}

\begin{document}

\begin{center}

{\bf Uniform Limit Theorem and Tail}\\

\vspace{3mm}

{\bf Estimates for parametric U \ - Statistics. }\\

\vspace{4mm}

{\bf E. Ostrovsky, L.Sirota.}\\

\vspace{5mm}

Department of Mathematics, Bar-Ilan University, Ramat-Gan, 59200, Israel.\\
e-mail: \ eugostrovsky@list.ru \\

\vspace{3mm}

Department of Mathematics, Bar-Ilan University, Ramat-Gan, 59200, Israel.\\
e-mail: \ sirota3@bezeqint.net \\

\vspace{4mm}

{\bf Abstract.} \par

\vspace{3mm}

\end{center}

 \  We deduce in this paper the sufficient conditions for weak convergence of centered and normed deviation of the $  \ U \ - \ $
statistics with values in the space of the real valued continuous function defined on some compact metric space. \par

\ We obtain also a non-asymptotic and  non-improvable up to multiplicative constant  moment
and exponential tail estimates  for  distribution for the uniform  norm of centered and naturally normed deviation of
$  U  \ - $ statistics by means of its martingale representation. \par

 \  Our results are formulated in a very popular and natural terms of metric entropy in the distance (distances)
generated by the introduced random processes (fields). \par

\vspace{4mm}

{\it Key words:} $  U  \ - $ statistics, kernel, rank, random variables (r.v.) and random fields (r.f.) (processes), distance,
compact metric space, space of continuous functions, uniform convergence, diameter, deviation, exponential Orlicz spaces,
weak distributions compactness and convergence, metric entropy and entropy integral, moment generating function, multiple stochastic integral,
martingales and  martingale differences, martingale representation,  Lebesgue-Riesz and Grand Lebesgue norm and spaces, symmetric function,
lower and upper moment and exponential estimates, moments, examples, natural functions, distance and norming; tails of distribution. \par

\vspace{3mm}


\vspace{4mm}

\section{ Introduction. Notations. Statement of problem.} \par

\vspace{3mm}

 \hspace{4mm}   Let $ (\Omega,F,{\bf P} ) $ be a probabilistic space, which will be presumed sufficiently
rich when we construct  examples (counterexamples),  $  T = \{ \ t \} $ be arbitrary compact (or semi-compact) metric space
relative {\it some } distance function   $ d = d(t,s), \ t,s \in T; $ the  concrete choice of the $  d(\cdot, \cdot)  $  will be
clarified below.\par

 \ Define $ I = I(n) = I(d; n) = \{i_1; i_2; \ldots; i_d \} $ the set of indices of the
form $ I(n) = I(d; n) = \{ \vec{i} \}  = \{i \} = \{i_1, i_2, \ldots, i_d \} $ such that $ 1 \le i_1 < i_2 <
i_3 < i_{d-1} < i_d \le n; \  J = J(n) = J(d; n) $ be the set of indices of
the form (subset of  $ I(d; n))  \  J(d; n) = J(n) =   \{ \vec{ j} \} = \{ j  \} = \{j_1; j_2; \ldots; j_{d-1} \} $
such that $ 1 \le j_1 < j_2 \ldots < j_{d-1} \le n- 1. $ \par

\vspace{3mm}

 \ Let also $ \{ \xi(i) \}, \ i = 1,2,\ldots, n, \ $ be independent identically distributed (i., i.d.)
random variables (r.v.) with values in the certain measurable space $  (X,S) , \  \Phi = \Phi(x(1), x(2), \ldots, x(d); t) $
be a symmetric measurable non-trivial numerical function (kernel) of $  d +1 $ variables:
$  \Phi: X^d \otimes T \to R, \ U(n,t) = U_n(t) = U(n, \Phi,d;t) =  $

$$
 U(n, \Phi,d; \{  \xi(i) \}; t ) = {n \choose d}^{-1} \sum_{I \in I(d,n)}
\Phi( \xi(i_1), \xi(i_2), \ldots, \xi(i(d)),t), \ n > d, \eqno(1.0)
$$
be a so-called $ U \ - $ statistic, which dependent also on the auxiliary  parameter $ t, \ t \in T. $ \par
 \ Denote $ \deg  \Phi = d, \ r = \rank  \Phi \in [ 1,2, \ldots, d-1 ], $ and suppose that $  r =\const, $ i.e. does not
dependent on the parameter $ t; $

$$
\Phi(t) := \Phi( \xi(1), \xi(2), \ldots, \xi(d); t), \hspace{4mm}
 \sigma(n) = \sigma_n(t) = \sqrt{ \Var (U_n(t))}, \eqno(1.1)
$$

\vspace{3mm}

  \ Let us introduce the following important notation. \par

\vspace{4mm}

{\bf Definition 1.1.} {\it Denote by}

$$
 \phi_n(t) \stackrel{def}{=} \frac{U_n(t) - {\bf E}U_n(t)}{ n^{r/2}} \eqno(1.2)
$$
{\it the centered and naturally normed deviation of our $  U \ -  $ statistics.} \par

\vspace{4mm}

 \ In the case when the rank of the considered $ \ U \ -  $ statistics is variable, one can understand in the definition (1.2)
as the capacity of the value $ r $  its {\it maximal} value:

$$
r := \max_{t \in T} \rank  \Phi(t).
$$

\vspace{4mm}

 \ Another possibility: to consider the behavior of the centered and normed values of our $ \ U \ -  $ statistics separately
 on the sets

$$
T_l = \{  t, \ \rank \Phi = l \}, \ l = 1,2, \ldots,r.
$$

 \ It is known that as $  n \to \infty  $ the finite-dimensional distributions of the r.f. $  \phi_n(\cdot) $ tends
to the  finite-dimensional distributions of the limited r.f. $ \phi(\cdot) = \phi_{ \infty }(\cdot), $  which is the  $ d \ -  $
multiple stochastic parametric integral, see  \cite{Korolyuk1}, chapter 4; \cite{Borovskikh1}.\par

\vspace{4mm}

 {\bf We will prove in this report under certain natural conditions the weak convergence in the space of all continuous functions
 $  C(T) $ the distribution the r.f. $ \ \phi_n \ $ to one for the r.f. $  \phi(\cdot), $ i.e. in the
 Prokhorov-Skorokhod sense. \par
 \ We  obtain also the non-asymptotical exponential decreasing, in general case,
  estimations for tail of the uniform norm for  $ \ \phi_n(\cdot). \ $  }

\vspace{4mm}

 \ On the other hands, we investigate the asymptotical and not asymptotical behavior of
 $ \ U \ - \ $ statistics in the (separable) Banach spaces $  C(T) $ of continuous
functions.  This problem may be named as $  UC \ - $ statistics, or equally $  UC \ - $ Limit Theorem  (UCLT),
 alike the notion of $ UH \ - $ statistics for Hilbert space valued
statistics and $ UB \ - $ statistics for arbitrary (separable) Banach space valued  statistics. \par
 \ The case of the Banach spaces with smooth norm and the Banach spaces of a finite type (co - type) $  p, \ p \ge 1 $
was considered in the books \cite{Korolyuk1}, \cite{Borovskikh1}. \par

 \ Note that the case $  d = r = 1 $ correspondent to the classical CLT in Banach space $  C(T),$  see \cite{Dudley1}, \cite{Ledoux1},
 \cite{Ostrovsky1} and reference therein. \par

\vspace{4mm}

 \ Recall that

$$
\sigma^2(n) = \Var(U(n)) \asymp n^{-r}, \ n \to \infty.
$$

\vspace{4mm}

 \ Here and in the future for any r.v. $ \ \eta \ $  the function $ \ T_{\eta} ( x) \ $ will be denote its {\it tail} function:

$$
 T_{\eta}(x) \stackrel{def}{=} \max( {\bf P}(\eta > x), {\bf P}(\eta < - x) ), \ x > 0.
$$

\vspace{4mm}

  \ The so-called {\it martingale representation } for the $ U \ - $ statistics as well as the
exact value for its variance $ \sigma^2(n) = \Var(U(n)) $ may be found, e.g. in \cite{Hoefding1}, \cite{Korolyuk1}, chapter 1.
It is very useful for the investigation of properties of these statistics.  \par

\vspace{3mm}

 \ We denote as usually the $ L(p) $ norm of the r.v. $ \eta $ as follows:

 $$
 |\eta|_p = \left[ {\bf E} |\eta|^p  \right]^{1/p}, \ p \ge 1; \eqno(1.3)
 $$

\vspace{4mm}

  \ Evidently, these weak convergence  may be applied for building  of an asymptotical confidence interval for unknown parameter
 by using the $ U \ -  $ statistics in the statistical estimation.  Namely, for any bounded and continuous functional
 $  F: C(T) \to R $

$$
\lim_{n \to \infty} {\bf E} F( \phi_n(\cdot) ) =  {\bf E} F( \phi(\cdot) ). \eqno(1.4)
$$
 \  In particular,

$$
\lim_{n \to \infty} {\bf P} \left( \max_{t \in T} | \ \phi_n(t) \ | > u \right) =
{\bf P} \left( \max_{t \in T} | \ \phi_{\infty}(t) \ | > u  \right), \ u > 0; \eqno(1.5)
$$

 \ Therefore, by the practical using

$$
 {\bf P} \left( \max_{t \in T} | \ \phi_n(t) \ | > u \right) \approx
{\bf P} \left( \max_{t \in T} | \ \phi_{\infty}(t) \ | > u  \right), \ u >> 1. \eqno(1.5a)
$$

 \ The exact asymptotical behavior as $ u \to \infty $ as well as non - asymptotical estimates for the right - hand side of
the relation (1.5a) is known, see  \cite{Piterbarg1}.\par

 \vspace{3mm}

 \ There are many famous works about this problem; the  next list is far from being complete:
  \cite{Borisov1}, \cite{Borovskikh1}, \cite{Pena1}, \cite{Gine1}, \cite{Hoefding1}, \cite{Korolyuk1}, \cite{Kwapien1},
  \cite{Ostrovsky10} etc.; see also the reference therein. \par
 \ Notice that in the  classical book  \cite{Korolyuk1}, see also \cite{Borovskikh1}, there are many examples of applying of
the theory of $ U \ -  $ statistics. A (relatively) new application, namely, in the modern adaptive estimation in the non - parametrical
statistics may be found in the article \cite{Bobrov1} and in the book \cite{Ostrovsky1}, chapter 5, section 5.13. \par

 \vspace{4mm}

\section{ Grand Lebesgue Spaces and  the exponential Orlicz spaces  for the analysis of random processes (fields). }

\vspace{4mm}

 \ Let $  \psi = \psi(p), \ p \in [2,b), \ b = \const, \ 2 < b \le \infty   $  (or $ p \in [2,b] ) $  be certain bounded from below:
 $ \inf \psi(p) \ge c = \const > 0 $  continuous
inside the {\it  semi - open  } interval $ [2, b) $ numerical function. We can and will suppose

$$
b = \sup \{p, \ \psi(p) < \infty \}, \eqno(2.0)
$$
so that $ \supp \psi = [2,b)  $ or $ \supp \psi = [2,b].$ The set of all such a functions will be denoted by
$ \Psi(b); \ \Psi:= \Psi(\infty). $ \par
 \ For each such a function $  \psi \in \Psi(b) $ we define

$$
\psi_d(p) \stackrel{def}{=}  \left[ \frac{p}{ \ln p} \right]^d \cdot \psi(p). \eqno(2.1)
$$
 \ Evidently,  $ \psi_d(\cdot) \in \Psi(b). $ \par
 \ This function  $ p \to (p/\ln p)^d  $
 play a very important role in the theory of summing of multiple martingale differences,  see
 \cite{Ostrovsky10}, \cite{Ostrovsky11}; for instance, for the estimation of moments for  sums of independent random variables, when
 $  d = 1; $ the correspondent estimate $ \psi_1(p) $ may be discovered by  H.P.Rosenthal \cite{Rosenthal1}.\par
 \ The multidimensional $ \ d \ge 2 \ $ martingale version of these  moment estimates was considered  in the authors preprint
\cite{Ostrovsky10}, which is based in turn on the  remarkable work belonging to  A.Osekowski \cite{Osekowski1}.\par

 \ By definition, the (Banach) space $  G\psi = G\psi(b)  $  consists on all the numerical valued random variables $ \{ \zeta \}  $
defined on our  probability space $ (\Omega,F,{\bf P} ) $ and having a finite norm

$$
||\zeta|| G\psi \stackrel{def}{=} \sup_{p \in (2,b)} \left[ \ \frac{|\zeta|_p}{\psi(p)} \ \right] < \infty. \eqno(2.2)
$$

 \ These spaces are suitable in particular for an investigation  of the random variables and the random processes (fields) with
exponential decreasing tails of distributions, the Central Limit Theorem in separable Banach spaces, study of Partial Differential
Equations etc., see e.g. \cite{Kahane1}, \cite{Buldygin1}, \cite{Kozatchenko1}, \cite{Ostrovsky1}, chapter 1,
\cite{Fiorenza1} - \cite{Fiorenza3}, \cite{Iwaniec1}  -  \cite{Iwaniec2} etc. \par
 \ More detail, suppose $  0 < ||\zeta|| := ||\zeta||G\psi < \infty. $ Define the function

$$
\nu(p) = \nu_{\psi}(p) \stackrel{def}{=} p \ln \psi(p), \ 2 \le p < b
$$
and put formally $ \nu(p) := \infty, \ p < 2  $ or $ p > b. $ Recall that the Young - Fenchel, or Legendre transform $ f^*(y) $
for arbitrary function $  f: R \to R $ is defined (in the one-dimensional case) as follows

$$
f^*(y) \stackrel{def}{=} \sup_x (x y - f(x)).
$$

 \ It is known that

 $$
 T_{\zeta}(y) \le \exp \left( - \nu_{\psi}^*(\ln (y/||\zeta||) )  \right), \ y > e \cdot ||\zeta||. \eqno(2.3)
 $$
 \ Conversely, if (2.3) there holds in the following version:

$$
T_{\zeta}(y) \le \exp \left( - \nu_{\psi}^*(\ln (y/K) )  \right), \ y > e \cdot K, \ K = \const > 0, \eqno(2.4)
$$
and the function $ \nu_{\zeta}(p), \ 2 \le p < \infty $ is positive, continuous, convex and such that

$$
\lim_{p \to \infty} \psi(p) = \infty,
$$
then $ \zeta \in G\psi $ and besides

$$
||\zeta||G\psi \le C(\psi) \cdot K. \eqno(2.5)
$$

\ Moreover, let us introduce the {\it exponential } Orlicz space $ L^{(M)} $ over the source probability space
$ (\Omega,F,{\bf P} ) $ with proper Young-Orlicz function

$$
M(u) := \exp \left(  \nu_{\psi}^*(\ln |u| )  \right), \ |u| > e
$$
or correspondingly

$$
M_d(u) := \exp \left(  \nu_{\psi_d}^*(\ln |u| )  \right), \ |u| > e
$$
and as ordinary $  M(u) = M_d(u) = \exp(C \ u^2) - 1, \ |u| \le e. $ It is known \ \cite{Ostrovsky11} that the $ G\psi $ norm
of arbitrary r.v. $  \zeta  $ is complete equivalent to the its norm in Orlicz space $ L^{(M)}: $

$$
||\zeta||G\psi \le C_1 ||\zeta||L^{(M)} \le C_2||\zeta||G\psi, \ 1 \le C_1 \le C_2 < \infty;
$$

$$
||\zeta||G\psi_d \le C_3 ||\zeta||L^{(M_d)} \le C_4||\zeta||G\psi_d, \ 1 \le C_3 \le C_4 < \infty.
$$

\vspace{4mm}

{\bf Example 2.1.} The estimate for the r.v. $ \xi $ of a form

$$
|\xi|_p  \le C_1 \ p^{1/m} \ \ln^r p, \ p \ge 2,
$$
 for some $ C_1 = \const < \infty, \ m = \const > 0, \ r = \const, \ $ is quite equivalent to the following tail estimate

$$
\exists   C_2 = C_2(C_1,m,r), \
T_{\xi}(x) \le  \exp \left\{ - C_2(C_1,m,r) \ x^m \ \ln^{-m r}x \right\}, \ x > e.
$$

\vspace{3mm}

 \ It is important to note that the inequality (2.4) may be applied still when the r.v. $ \xi $ does not have the
exponential  moment, i.e. does not satisfy the famous Kramer's condition. Namely, let us consider the next example. \\

\vspace{3mm}

 {\bf Example 2.2.} Define the following $ \Psi \ -  $ function.

$$
\psi_{[\beta]}(p) := \exp \left( C_3 \ p^{\beta} \right), \  p \in [2, \infty), \ \beta = \const > 0.
$$

 \ The r.v. $  \xi $ belongs to the space $ G \psi_{[\beta]} $ if and only if

$$
T_{\xi}(x) \le \exp \left( - C_4(C_3,\beta) \ [\ln (1 + x)]^{1 + 1/\beta}   \right), \ x \ge 0.
$$

 \ See also \cite{Ostrovsky4}. \par

\vspace{3mm}

 \ Let us return to the source problem. Assume that there exists certain function $  \psi(\cdot) \in \Psi(b), \   b = \const \in (2, \infty)  $
such that $ \Phi(t) \in G\psi(b) $ uniformly in $  t: $

$$
\sup_{t \in T} || \ \Phi(t) \ ||G\psi < \infty.  \eqno(2.6)
$$

 \ For instance, this function may be picked by the following {\it natural} way:

$$
\psi_{\Phi}(p) := \sup_{t \in T} |\Phi(t)|_p, \eqno(2.7)
$$
if of course there exists and is finite at last for some value $  p  $ greatest than 2,  obviously,
with the appropriate value $  b. $\par

 \ Evidently,

$$
\sup_{t \in T} || \ \Phi(t) \ || G\psi_{\Phi} = 1.
$$

 \ Assume now that {\it some } separable  stochastic continuous random  field  $ \ \eta(t) \ $ satisfies the condition (2.6),
 relative some $ \psi(\cdot) \ -  $ function, $ \psi(\cdot)  \in \Psi(b), \ $ i.e.

$$
\sup_{t \in T} || \ \eta(t) \ ||G\psi < \infty.  \eqno(2.8)
$$

 \ Define as ordinary $ d = d(s,t) = d_{\eta}(s,t)= d_{\psi, \eta}(s,t)  $ is
 the following natural bounded distance function (more precisely, semi - distance) on the set $  \ T: \  $

$$
d(t,s) = d_{\eta}(t,s) \stackrel{def}{=} || \ \eta(t) - \eta(s)  \  || G\psi, \eqno(2.9)
$$

and as ordinary

$$
\diam(T,d) \diam_d*(T) = \sup_{t,s \in T} d(t,s) < \infty.
$$

\vspace{3mm}

 \ Let us introduce for any subset $ V, \ V \subset T $ the so-called
{\it entropy } $ H(V, d, \epsilon) = H(V, \epsilon) $ as a logarithm
of a minimal quantity $ N(V,d, \epsilon) = N(V,\epsilon) = N $
of a balls $ S(V, t, \epsilon), \ t \in V: $
$$
S(V, t, \epsilon) \stackrel{def}{=} \{s, s \in V, \ d(s,t) \le \epsilon \},
$$
which cover the set $ V: $
$$
N = \min \{M: \exists \{t_i \}, i = 1,2,…, M, \ t_i \in V, \ V
\subset \cup_{i=1}^M S(V, t_i, \epsilon ) \}, \eqno(2.10)
$$
and we denote also
$$
H(V,d,\epsilon) = \log N; \ S(t_0,\epsilon) \stackrel{def}{=}
 S(T, t_0, \epsilon), \ H(d, \epsilon) \stackrel{def}{=} H(T,d,\epsilon). \eqno(2.11)
$$
 \ It follows from the famous Hausdorff's theorem that
$ \forall \epsilon > 0 \ \Rightarrow H(V,d,\epsilon)< \infty $ iff the
metric space $ (V, d) $ is precompact set, i.e. is the bounded set with
compact closure.\par

 \ Define the function

 $$
  v_{\psi}(x) \stackrel{def}{=} \inf_{p \in (2,b)} \left( \frac{x}{p} + \ln \psi(p) \right), \eqno(2.12)
 $$
and the following so-called {\it entropy} integral

$$
I = I(\psi,d) \stackrel{def}{=} \int_0^1 \exp \left( v_{\psi}(H(T,d,\epsilon))  \right) \ d \epsilon. \eqno(2.13)
$$

\ It is proven in \cite{Ostrovsky1},  pp. 172-176 that if $   I(\psi,d) < \infty, $ then almost all trajectories  of the r.f.
$ \ \eta(\cdot) \ $ are  $  d \ -  $ continuous.\par
 \ Furthermore,

$$
 || \ | \sup_{t \in T} \ | \eta(t) \ | \ ||G\psi \le C(H(\cdot,d(\cdot),\cdot), \diam_d(T), I(\psi,d)),
\eqno(2.13a)
$$
which allows in turn to obtain in general case the exponential bounds for the distribution of
$ \ \sup_{t \in T} | \ \eta(t) \ | $ in accordance with estimates (2.3) and (2.4). \par

\vspace{4mm}

 \ Moreover,  let $ \eta_n(t), \ t \in T, \ n = 1,2,\ldots $ be certain {\it sequence} of separable random fields such that
for some non-random point $ t_0 \in T $ the sequence of distributions of the one-dimensional r.v.  $ \eta_n(t_0) $ is
weak compact on the real line. Define the natural $ \  \psi \ - \ $ function for the family of r.v. $ \eta_n(t): $

$$
\theta(p) := \sup_n \sup_t | \ \eta_n(t) \ |_p,
$$
and suppose its finiteness for some set $ p \in (2,b). $ Introduce as before  the (bounded) distance function

$$
r(t,s) = r_{\theta}(t,s) \stackrel{def}{=} \sup_n || \ \eta_n(t) - \eta_n(s) \  ||G\theta.
$$

 \ If the following entropy integral converges

$$
J = J(\theta, r_{\theta}) \stackrel{def}{=} \int_0^1 \exp \left( v_{\theta}(H(T,r,\epsilon))  \right) \  d \epsilon < \infty,
\eqno(2.14)
$$
then the family of distributions generated in the space $ C(T,r) $ by the random fields $ \eta_n(\cdot)  $ is weakly compact. \par

 \ Herewith  as before

$$
\sup_n || \ \sup_{t \in T} |\eta_n(t)|  \ ||G\theta \le C \left\{r, \theta, J,  \diam_{\theta}(T) \right\} < \infty, \eqno(2.14a)
$$
 which allows us to obtain in general case the {\it exponential decreasing} tail estimate for the r.v.
 $ \zeta(n):=  \sup_{t \in T} |\eta_n(t)|:  $

 $$
 \sup_n T_{\zeta(n)}(y) \le \exp \left\{- \nu_{\theta}^*\left[\ln \left(y/C_2(r,\theta,J, \diam_{\theta}(T)) \right) \right] \right\}, \ y > C_2 e.
 $$

\vspace{4mm}

\ {\bf Example 2.3.} Let $  \eta = \eta(t), \ t \in T $ be separable random field such that for some $  \ b = \const \ge 1 \  $

$$
\sup_t | \ \eta(t) \ |_b < \infty.
$$
 \ Define the following bounded natural distance $  d_{(b)} =  d_{(b)}(t,s)  $ on the set $ T: $

$$
d_{(b)}(t,s) := | \ \eta(t) - \eta(s) \ |_b.
$$

 \ The integral (2.13) has a form

$$
I = C \int_0^1 N^{1/b}(T,d_{(b)},\epsilon) \ d \epsilon. \eqno(2.15)
$$
 \ The finiteness of the last integral guarantee the $ d_{(b)} (\cdot) $ continuity of the r.f. $  \eta(\cdot) $ with probability one:

$$
{\bf P} \left(\eta(\cdot) \in C(T,d_{(b)}) \right) = 1
$$
and herewith

$$
 | \ \sup_{t \in T} \ | \eta(t) \ | \ |_b \le C_1(N(T,d_b(\cdot),\cdot), \diam(d_b(T), I(\psi,d)). \eqno(2.15a)
$$

 \ We obtained the remarkable Pisier's condition, see \cite{Pizier1}, \cite{Pizier2}. \par

 \ Let in addition the set $ T  $ be a closure of a non-empty convex bounded set in the whole space $  R^k, \ k = 1,2,\ldots.:
 \ T \subset R^k. $ Assume that the distance $  d_b(t,s) $ is such that

$$
d_b(t,s) \le C_1 \ |t - s|^{\alpha}, \ t,s \in T, \ \alpha = \const \in (0,1],
$$
where $ |t| $ denotes the ordinary Euclidean norm. Then

$$
N(T,d_b,\epsilon) \le C_2 \ \epsilon^{k/\alpha}, \ \epsilon \in (0,1).
$$

 \ The condition (2.14) is satisfied iff $ \alpha > k/b.  $ More detail:

$$
\exists t_0 \in T  \ \Rightarrow \sup_n  {\bf E} | \ \eta_n(t_0) \ |^b < \infty
$$
and
$$
\exists \beta = \const > k, \ C = \const < \infty, \ \Rightarrow
\sup_n {\bf E} | \ \eta_n(t) - \eta_n(s) \ |^b \le C \ |t - s|^{\beta}.
$$

 \ We obtained in fact the noticeable Kolmogorov-Slutsky condition. \par

 \vspace{4mm}

\ The  {\it  minimal} value of the variable $ \alpha/k $ is said to be {\it  entropy dimension } of the space $ T $
relative the distance $ d = d_b, $  write

$$
\dim_{d}T := \inf (k/\alpha).
$$

\vspace{4mm}

\ {\bf Example 2.4.} Let $  \eta = \eta(t), \ t \in T $ be separable random field belonging uniformly relative the parameter $ t \in T $
to the space $ G \psi_{m,r} $ of the form

$$
\psi_{m,r}(p) \stackrel{def}{=} p^{1/m} \ \ln^r p, \ p \ge 2, \ m = \const > 0, \ r = \const \in R; \eqno(2.16)
$$
here $  b = \infty, $ uniformly in $ t, \ t \in T. $ This implies by definition

$$
\sup_{t \in T} | \ \eta(t) \ |_p < C  p^{1/m} \ \ln^r p = C \psi_{m,r}(p), \ p \ge 2, \eqno(2.17)
$$
or

$$
\sup_{t \in T} || \ \eta(t) \ ||G\psi_{m,r} < \infty. \eqno(2.17a)
$$

 \ Define as before  the bounded natural distance $ \rho_{m,r}  $ on the set $  T  $

$$
\rho_{m,r}(t,s) := || \ \eta(t) - \eta(s) \  ||G\psi_{m,r}. \eqno(2.18)
$$

 \ The condition (2.14) takes the form

$$
J(m,r):= \int_0^1 H^{1/m}(T, \rho_{m,r}, \epsilon) \ | \ \ln H(T, \rho_{m,r}, \epsilon) \ |^r \ d \epsilon < \infty  \eqno(2.19)
$$
and guarantee us the $ \rho_{m,r}(\cdot, \cdot)  $ continuity of the r.f. $ \eta(\cdot) $ almost surely. Moreover,

$$
{\bf P} \left( \sup_{t \in T} |\eta(t)| > u \right) \le
$$

$$
\exp \left\{ \ - C_2 \left[m,r,J(m,r), \diam_{\rho_{m,r}}(T)\right] \ u^m \ \ln^{-mr}u \ \right\}, \ u \ge e. \eqno(2.20)
$$

 \ If for instance the r.f. $ \eta(t) $ is centered Gaussian, then we can choose $ m = 2 $ and $  r = 0. $ The condition (2.19)
takes the form

$$
J^o:= \int_0^1 H^{1/2}(T, \rho_{2,0}, \epsilon) \ d \epsilon < \infty, \eqno(2.21)
$$
see \cite{Fernique1},  \cite{Dudley2}.\par

 \ If in addition the r.f.  $ \eta(t) $ is stationary and $ T = [0,1] $ or $ T = [0, 2 \pi], $ then this condition is not only
 sufficient, but it is necessary still for the boundedness a.e. of the r.f. $ \eta(\cdot), $ see \cite{Fernique1}. \par

\vspace{4mm}

\ {\bf Example 2.5.} Let $  \eta = \eta(t), \ t \in T $ be separable random field belonging uniformly in $ t \in T $ to the space
$ G \psi_{[\beta]}, $ see example (2.2). Moreover, we suppose


$$
\sup_{t \in T} T_{\eta(t)}(x) \le
\exp \left( - C \ [\ln (1 + x)]^{1 + 1/\beta}   \right), \ x \ge 0,\eqno(2.22)
$$
where $  0 < C = \const < \infty. $  The natural distance $ q(t,s)  $ may be defined as ordinary

$$
q(t,s) := || \ \eta(t) - \eta(s)  \ ||G\psi_{[\beta]}. \eqno(2.23)
$$

 \ The condition (2.13) takes the form

$$
L:= \int_0^1 \exp \left( H^{\beta/(\beta + 1)} (T,q,\epsilon)  \right) \ d \epsilon < \infty. \eqno(2.24)
$$

 \ Furthermore, if this condition is satisfied, then for all the positive values $  x  $

$$
{\bf P} \left( \sup_{t \in T} |\ \eta(t) \ | > x \right) \le \exp \left( - C_3(\beta,L, \diam_q(T)) \ [\ln (1 + x)]^{1 + 1/\beta} \right).
\eqno(2.25)
$$

\vspace{4mm}

\section{Main result. }

\vspace{4mm}

 \hspace{4mm} Let us return to formulated above in first section problem. Suppose that for the r.f. $ \Phi^o(t)= \Phi(t) - {\bf E}\Phi(t) $
there exists certain $ \psi \ -  $ function from the set $  \Psi(b), \ 2 < b \le \infty $  such that

$$
\sup_{t \in T} || \ \Phi^o(t)  \  ||G\psi < \infty, \eqno(3.1)
$$
or equally

$$
\sup_{t \in T} | \ \Phi^o(t)  \  |_p < C \ \psi(p), \ p \in [2,b). \eqno(3.1a)
$$

 \ For instance, the function $ \psi(\cdot) $ may be picked as a natural function for the family $  \{ \Phi^o(t) \}, \ t \in T, $
if of course there exists for some value $ b > 2: $

$$
\tilde{\psi}(p) := \sup_{\in T} | \ \Phi^o(t) \ |_p, \ 2 \le p < b. \eqno(3.1a)
$$

\vspace{4mm}

 \ {\it We can and will agree  in what follows  without loss of generality} $ {\bf E} \Phi(t) = 0, $ so that $ \Phi^o(t) = \Phi(t).  $ \par

\vspace{4mm}

 \ Let us introduce the following bounded natural semi-distance $  w(t,s) = w_{\Phi}(t,s):  $

$$
w(t,s) = w_{\Phi}(t,s) \stackrel{def}{=} || \ \Phi(t) - \Phi(s)  \ ||G\psi. \eqno(3.2)
$$

 \ Define a new $ \psi \ - $ function $  \  \tau = \tau(p), \ 2 \le p < b: $

$$
\tau(p) = \tau_{\Phi,d}(p) := \left[\frac{p}{\ln p}  \right]^d \ \psi(p). \eqno(3.3)
$$

\ It is proven in particular in \cite{Ostrovsky11} that

$$
\sup_n | \ \phi_n(t) \ |_p \le C \ \left[\frac{p}{\ln p}  \right]^d \cdot | \ \Phi \ |_p \le
$$

$$
C \ \left[\frac{p}{\ln p}  \right]^d \cdot\psi(p) = C \ \tau(p), \eqno(3.4)
$$
or equally

$$
\sup_n \sup_{t \in T} || \ \phi_n(t) \ ||G\tau \le C < \infty.\eqno(3.4a)
$$

 \ Since the sequence of differences $ \Delta U_n = U_n(t) - U_n(s) $ forms also the $ U \ - $ statistics, we conclude analogously
to the inequality (3.4a)

$$
\sup_n || \ \phi_n(t) - \phi_n(s) \ ||G\tau \le C \ w(t,s).\eqno(3.5)
$$

  \ It remains to apply the proposition (2.14) to obtain the following result.\\

\vspace{4mm}

{\bf Theorem 3.1.} {\it Suppose in addition}

$$
K = K(w,\tau) := \int_0^1 \exp \left\{ v_{\tau}(H(T,w,\epsilon)) \right\} \ d \epsilon < \infty. \eqno(3.6)
$$
 \ {\it Then then the family of distributions generated in the space $ C(T,w) $ by the $ w \ - $ continuous with probability one random fields
$ \phi_n(\cdot)  $ is weakly compact, and in addition  as before}

$$
\sup_n || \ \sup_{t \in T} |\phi_n(t)|  \ ||G\tau \le C(w, \theta, K, \diam_w(T)) < \infty, \eqno(3.7)
$$
{\it with the correspondent tail estimate}

 $$
\sup_n T_{\sup_t |\phi_n(t)|}(y) \le
\exp \left\{ - \nu_{\tau}^* \left[ \ \ln (y/C_2(w,\theta,K, \diam_w(T)) \ ) \ \right] \ \right\}, \ y > C_2 e. \eqno(3.7a)
 $$

\vspace{4mm}

{\bf Example 3.1.} Assume that the (random) function $ \Phi(\cdot) $ allows the following estimate

$$
\sup_{t \in T} | \ \Phi(t)  \  |_b  < \infty
$$
for some constant $  b = \const \ge 2. $ The distance function $  w = w(t,s) $ may be introduced as above

$$
 w(t,s) = |  \ \Phi(t) -\Phi(s)  \  |_b,
$$
 and the condition (3.6) coincides with the Pisier's condition

$$
\int_0^1 N^{1/b} (T,w,\epsilon) \ d \epsilon < \infty.
$$
 \ If this condition there holds, we deduce as a consequence:

$$
\sup_n | \ \sup_{t \in T} |\phi_n(t)|  \ |_b < \infty.
$$

\vspace{4mm}

{\bf Example 3.2.} Impose on the kernel-function $ \Phi(t)  $ the following condition:

$$
 \sup_{t \in T} || \ \Phi(t) \ || G \psi_{m,r}  < \infty, \ m = \const > 0, \ r = \const \in R, \eqno(3.8)
$$
or equally

$$
 \sup_{t \in T} | \ \Phi(t) \ |_p \le C p^{1/m} \ \ln^r p, \ m = \const > 0, \ r = \const \in R, \eqno(3.8a)
$$

 \ The correspondent (bounded) natural distance $ z_{\Phi} =  z_{\Phi}(t,s)   $ function has the form

$$
z_{\Phi}(t,s) := || \  \Phi(t) - \Phi(s) \ ||G \psi_{m,r}.
$$

\ As we knew,  see  \cite{Ostrovsky11},

$$
\sup_n | \ \phi_n(t) \ |_p \le C \ \left[\frac{p}{\ln p}  \right]^d \cdot | \ \Phi \ |_p \le
$$

$$
C \ \left[\frac{p}{\ln p}  \right]^d \cdot p^{1/m} \ \ln^r p = C \ p^{(1 + dm)/m} \ \ln^{r - d}p, \ p \ge 2,
$$
or equally

$$
\sup_n \sup_{t \in T} || \ \phi_n(t) \ ||G\psi_{l, g}  \le C < \infty, \eqno(3.9)
$$
where

$$
l= l(m,d) = m/(1 + d m), \ g = g(r,d) = r - d.
$$

 \ The condition (3.6) takes the form

$$
K_{m,r,\Phi} := \int_0^1 H^{ (1 + d m)/m}(T, z_{\Phi}, \epsilon) \ | \ \ln H( T, z_{\Phi}, \epsilon ) \ |^{r - d} d \epsilon < \infty.
\eqno(3.10)
$$

 \ If the condition (3.10) there holds, we conclude by virtue of theorem 3.1  that
 the family of distributions generated in the space $ C(T,z_{\phi}) $ by the random fields $ \phi_n(\cdot)  $ is weakly compact, and
in addition

$$
\sup_n || \ \sup_{t \in T} |\phi_n(t)|  \ ||G\psi_{l(m,d), g(r,d)} < \infty \eqno(3.11)
$$
with correspondent non-asymptotical tail estimate (2.3), (2.4) for the uniform norm  $ \ \sup_{t \in T} |\phi_n(t)|: $ \par

$$
\sup_n  {\bf P} ( \sup_{t \in T} |\phi_n(t)|  > u)  \le \exp \left( - C_2 u^{l(m,d)} \ (\ln u)^{-l(m,d) \ g(r,d)} \right),
 \ u > e.\eqno(3.12)
$$

 \ Let for example again the set $ T  $ be a closure of a non-empty convex bounded set in the whole space $  R^k, \ k = 1,2,\ldots.:
 \ T \subset R^k. $ Assume that the distance $  z_{\Phi}(t,s) $ is such that

$$
z_{\Phi}(t,s) \le C_1 \ |t - s|^{\alpha}, \ t,s \in T, \ \alpha = \const \in (0,1],
$$
where $ |t| $ denotes the ordinary Euclidean norm. Then

$$
N(T,d_b,\epsilon) \le C_2 \ \epsilon^{\alpha/k}, \ \epsilon \in (0,1).
$$
 \ The condition (3.6) is satisfied for arbitrary value $ \alpha > 0.  $ \par

\vspace{4mm}

\ {\bf Example 3.3.} Let the continuous kernel-function $  \Phi = \Phi(t) $ be such that
$  \Phi = \Phi(t), \ t \in T $  belong uniformly in $ t \in T $ to the space
$ G \psi_{[\beta]}, \ \beta = \const > 0; $ see examples (2.2) and (2.5). Moreover, we suppose

$$
\exp \left( - C_1 \ [\ln (1 + x)]^{1 + 1/\beta}   \right) \le
\inf_{t \in T} T_{\Phi(t)}(x) \le \sup_{t \in T} T_{\Phi(t)}(x) \le
$$

$$
\exp \left( - C_2 \ [\ln (1 + x)]^{1 + 1/\beta}   \right), \ x \ge 0, \eqno(3.13)
$$
where $  0 < C_2 = \const \le C_1 = \const < \infty. $   In particular,

$$
 \sup_{t \in T} | \ \Phi(t) \ |_p \le \exp \left(  C_3 p^{\beta} \right). \eqno(3.13a)
$$

 \ The natural (and bounded) distance $ q_{\Phi,\beta}(t,s)  $ may be defined as ordinary

$$
q_{\Phi,\beta}(t,s) := || \ \Phi(t) - \Phi(s)  \ ||G\psi_{[\beta]}. \eqno(3.14)
$$

\ As we knew,  see  \cite{Ostrovsky11},

$$
\sup_n | \ \phi_n(t) \ |_p \le C \ \left[\frac{p}{\ln p}  \right]^d \cdot | \ \Phi \ |_p \le
$$

$$
C \ \left[\frac{p}{\ln p}  \right]^d \cdot \exp \left(  C_3 \ p^{\beta} \right) \le \exp \left(  C_4(\beta) \ p^{\beta} \right), \ p \ge 2,
$$
or equally

$$
\sup_n \sup_{t \in T} || \ \phi_n(t) \ ||G\psi_{[\beta]}  \le C_5 < \infty.
$$

 \ Analogously,

$$
\sup_n  || \ \phi_n(t) - \phi_n(s) \ ||G\psi_{[\beta]}  \le C_5 \ q_{\Phi,\beta}(t,s).
$$

 \ The condition (3.6) takes the form

$$
L = L(\Phi,\beta):= \int_0^1 \exp \left( H^{\beta/(\beta + 1)} (T,q_{\Phi,\beta},\epsilon)  \right) \ d \epsilon < \infty. \eqno(3.15)
$$

 \ If the last condition there holds, we conclude by virtue of theorem 3.1  that
 the family of distributions generated in the space $ C(T,d_{\Phi,\beta}) $ by the random fields $ \phi_n(\cdot)  $ is weakly compact, and
in addition

$$
\sup_n || \ \sup_{t \in T} |\phi_n(t)|  \ ||G\psi_{[\beta]} = C_6(\Phi,\beta,L, \diam \left(q_{\Phi,\beta},T \right))  < \infty \eqno(3.16)
$$
with the correspondent non-asymptotical tail estimate (2.3), (2.4) for the uniform norm  $ \ \sup_{t \in T} |\phi_n(t)|: $ \par

$$
\sup_n  {\bf P} ( \sup_{t \in T} |\phi_n(t)|  > u)  \le
$$

$$
\exp \left( - C_7 (\Phi,\beta,L, \diam_{q_{\Phi,\beta}}, \ T) \ (\ln(1 + u))^{1 + \beta} \right), \ u > 0. \eqno(3.17)
$$

\vspace{4mm}

 \ Notice that the condition (3.15) is satisfied if for instance

$$
\dim_{q_{\Phi,\beta}} T < \infty.
$$

\vspace{4mm}

 \ Let us deduce the {\it lower bound}  for this tail probability. We conclude taking into account the condition (3.13)  and
choosing arbitrary (non-random) value $  t_o \in T  $

$$
\sup_n  {\bf P} ( \sup_{t \in T} |\phi_n(t)|  > u)  \ge  {\bf P} ( |\phi_1(t_o)|  > u)  =
$$

$$
{\bf P} ( |\Phi(t_o)|  > u) \ge
\exp \left( - C_1  \ (\ln(1 + u))^{1 + \beta} \right), \ u > 0.  \eqno(3.18)
$$

 \ Wherein the set $  T $ may consists on the unique point $ t_o; $ the condition (3.15) is not necessary for our lower bound. \par

 \ Thus, the estimation  (3.7), (3.7a) of theorem 3.1 is essentially, i.e. up to multiplicative constant, not improvable. \par

\vspace{4mm}

\section{Concluding remarks. }

\vspace{4mm}

 \hspace{4mm} {\bf A.}  CLT in the space of continuous functions as a particular case. \par
 \ Notice that the case $  d = r = 1 $  correspondent to the classical Central Limit Theorem in the space of continuous functions,
see  \cite{Dudley1}, \cite{Ostrovsky1}. Authors hope that our result, indeed, theorem 3.1, is some extension of the propositions
obtained therein.\par

 \vspace{3mm}

 {\bf B.} It is interest, by our opinion, to obtain  analogous estimates for dependent source random variables,
for instance, for martingales or mixingales. Some preliminary results in this directions may be found in \cite{Borisov1},
 \cite{Gine1}. \\

 \vspace{3mm}

 {\bf C.} We do not aim to derive the best possible  values of appeared in this report constants.
It remains to be done. \\

\vspace{3mm}

{\bf D.} Obviously, the case of the so-called  $ V \ - $ statistics  may be investigated analogously. \\

\vspace{3mm}

{\bf E.} Perhaps, more general results in the considered in this report problem may be obtained by means of more modern
technic, namely, through the so-called {\it majorizing measure} method, see. e.g.
\cite{Fernique1}, \cite{Ledoux1}, \cite{Talagrand1}-\cite{Talagrand4}.\par

\end{document}